\documentclass[12pt]{article}
\usepackage{latexsym}
\usepackage{amssymb}
\usepackage{amsmath}

\begin{document}

\newtheorem{theorem}{Theorem}
\newtheorem{lemma}[theorem]{Lemma}
\newtheorem{corollary}[theorem]{Corollary}
\newcommand{\forces}{\ \rule{.12mm}{2.9mm}\!\vdash}
\newcommand{\bnep}{\dot{\varepsilon}}
\newcommand{\vs}{\vspace{12pt}}
\newcommand{\A}{{\mathcal A}}
\newcommand{\B}{{\mathcal B}}
\newcommand{\nof}{{\mathbb P}}
\newcommand{\gnof}{{\mathbb P} \otimes \dot{\mathbb Q}}
\newcommand{\Q}{{\mathbb Q}}
\newcommand{\R}{{\mathbb R}}
\newcommand{\nofS}{{\mathbb S}}
\newcommand{\Z}{{\mathbb Z}}
\newcommand{\cons}{\textbf{C}}
\newcommand{\Aut}{\mathrm{Aut}}
\newcommand{\con}{\mathrm{con}}
\newcommand{\dom}{\mathrm{dom}}
\newcommand{\neat}{\mathcal{N}}
\newcommand{\proof}{\textbf{\noindent Proof}}

\title{Naturality and definability II}
\author{Wilfrid Hodges and Saharon Shelah}
\date{Draft, 13 November 2000}
\maketitle

In two papers \cite{ho:1} and \cite{hosh:1} we noted that in common practice
many algebraic constructions are defined only `up to isomorphism' rather
than explicitly.  We mentioned some questions raised by this fact, and we
gave some partial answers. The present paper provides much fuller answers,
though some questions remain open.  Our main result, Theorem~\ref{th:4},
says that there is a transitive model of Zermelo-Fraenkel set theory with
choice (ZFC) in which every explicitly definable construction is `weakly
natural' (a weakening of the notion of a natural transformation).  A
corollary is that there are models of ZFC in which some well-known
constructions, such as algebraic closure of fields, are not explicitly
definable.  We also show (Theorem~\ref{th:2}) that there is no transitive
model of ZFC in which the explicitly definable constructions are precisely
the natural ones.

Most of this work was done when the second author visited the first at Queen
Mary, London University under SERC Visiting Fellowship grant GR/E9/639 in
summer 1989, and later when the two authors took part in the Mathematical
Logic year at the Mittag-Leffler Institute in Djursholm in September 2000.
The second author proposed the approach of section \ref{se:3} on the first
occasion and the idea behind the proof of Theorem \ref{th:4} on the second.
Between 1975 and 2000 the authors (separately or together) had given some
six or seven false proofs of versions of Theorem \ref{th:4} or its negation.

The authors thank Ian Hodkinson for his invaluable help (while research
assistant to Hodges under SERC grant GR/D/33298) in unpicking some of the
earlier false proofs.

\section{Constructions up to isomorphism} \label{se:1}

To make this paper self-contained, we repeat or paraphrase some definitions
from \cite{hosh:1}.

Let $M$ be a transitive model of ZFC (Zermelo-Fraenkel set
theory with choice). By a \emph{construction} (in $M$) we mean a triple
${\bf C}=\langle \phi_1,\phi_2,\phi_3\rangle$ where
\begin{enumerate}
\item $\phi_1(x)$, $\phi_2(x)$ and $\phi_3(x)$ are formulas of
the language of set theory, possibly with parameters from $M$;
\item $\phi_1$ and $\phi_2$ respectively define first-order
languages $L$ and $L^-$ in $M$;  every symbol of $L^-$ is a symbol
of $L$, and the symbols of $L \setminus L^-$ include a 1-ary relation
symbol $P$;
\item the class $\{a : M \models \phi_3(a)\}$ is in $M$ a class of
$L$-structures, called the \emph{graph} of \textbf{C};
\item if $B$ is in the graph of \textbf{C} then $P^B$, the set of
elements of $B$ satisfying $Px$, forms the domain of an $L^-$-structure
$B^-$ inside $B$;  the class of all such $B^-$ as $A$ ranges over
the graph of \textbf{C} is called the \emph{domain} of \textbf{C};
\item the domain of \textbf{C} is closed under isomorphism; and
if $A, B$ are in the graph of \textbf{C} then every isomorphism
from $A^-$ onto $B^-$ extends to an isomorphism from $A$ onto $B$.
\end{enumerate}

A typical example is the construction whose domain is the class of
fields, and the structures $B$ in the graph are the algebraic closures of
$B^-$, with $B^-$ picked out by the relation symbol $P$. The algebraic
closure of a field is determined only up to isomorphism over 
the field;  in the terminology below, algebraic closures are 
`representable' but not `uniformisable'.  (What we called
definable in \cite{hosh:1} we now call \emph{uniformalisable};
the new term is longer, but it is less misleading because it
agrees better with the common mathematical use of these words.)

We say that the construction $\cons$ is $X$-\emph{representable} (in $M$) if
$X$ is a set in $M$ and all the parameters of $\phi_1$, $\phi_2$, $\phi_3$
lie in $X$.  We say that $\cons$ is \emph{small} if the domain of $\cons$
(and hence also its graph) contains only a set of isomorphism types of
structures.

An important special case is where the domain of $\cons$ contains exactly
one isomorphism type of structure; in this case we say $\cons$ is
\emph{unitype}.

The map $B^- \mapsto B$ on the domain of a construction \textbf{C} is in
general not single-valued; but by clause (5) it is single-valued up to
isomorphism over $B^-$. We shall say that \textbf{C} is
\emph{uniformisable} if its graph can be uniformised, i.e.\ there is a
formula $\phi_4(x,y)$ of set theory (the \emph{uniformising formula}) such
that 
\begin{quote}
for each $A$ in the domain of \textbf{C} there is a unique $B$ such
that $M \models \phi_4(A,B)$, and this $B$ is an $L$-structure in
the graph of \textbf{C} with $A = B^-$.
\end{quote}
We say that \textbf{C} is $X$-\emph{uniformisable} (in $M$) if there is such
a $\phi_4$ whose parameters lie in the set $X$.

\section{Splitting, naturality and weak naturality}  \label{se:2}
Let $\nu : G \to H$ be a surjective group homomorphism.  A \emph{splitting}
of $\nu$ is a group homomorphism $s: H \to G$ such that $\nu s$ is the
identity on $H$.  We say that $\nu$ \emph{splits} if it has a splitting.

For our Theorem \ref{th:4} we shall need a weakening of these notions.  A
stronger version of Theorem \ref{th:4} would make this unnecessary, but we
do not know whether the stronger version is true.

Let $\nu: G \to H$ be as above.  By a \emph{weak splitting} of $\nu$ we mean
a map $s: H \to G$ such that 
\begin{enumerate}
\item[(a)] $\nu s$ is the identity on $H$;
\item[(b)] there is a commutative subgroup $G_0$ of $G$ 
such that if $f_1$, \ldots, $f_k$ are elements of $H$ for which
$f_1^{\varepsilon_1} \ldots f_k^{\varepsilon_k} = 1$ (where
$\varepsilon_i$ is each either 1 or $-1$),
then $s(f_1)^{\varepsilon_1} \ldots s(f_k)^{\varepsilon_k} \in G_0$.
\end{enumerate}
If we strengthened this definition by requiring $G_0$ to be $\{1\}$, it
would say exactly that $s$ is a splitting of $\nu$.  In particular every
splitting is a weak splitting.  We say that $\nu$ \emph{weakly splits} if it
has a weak splitting.

Suppose $s$ is a weak splitting of $\nu$.  Then there is a smallest group
$G_0$ as in (b); it is the group consisting of the words
$s(f_1)^{\varepsilon_1} \ldots s(f_k)^{\varepsilon_k}$ as in (b).  This
group $G_0$ has the property that if $g$ is in $G_0$ and $f$ is in $H$ then
$s(f)^{-1}gs(f)$ is also in $G_0$.  So the normaliser of $G_0$ in $G$
contains the image of $s$.
\vs

\textbf{Example 1}. Let $G$ be the multiplicative group of $3\times 3$ upper
unitriangular matrices over the ring $\Z /(8\Z )$.  Let $H$ be the
corresponding group over $\Z /(2\Z )$, and let $\nu: G \to H$ be the
canonical surjection.  We show that $\nu$ doesn't weakly split.

Suppose for contradiction that $s$ is a weak splitting of $\nu$.  Let
$g_1, g_2$ be the two matrices
\[ 
g_1 = \left( \begin{array}{ccc}
1 & 1 & 0  \\
0 & 1 & 0  \\
0 & 0 & 1    \end{array} \right),  \ \   
g_2 = \left( \begin{array}{cccc}
1 & 0 & 0  \\
0 & 1 & 1  \\
0 & 0 & 1   \end{array} \right)
\]
in $G$, and write $f_1 = \nu(g_1)$, $f_2 = \nu(g_2)$.  Now $f_1^2 = f_2^2 =
1$ in $H$, so the weak splitting property tells us that $s(f_1)^2$ and
$s(f_2)^2$ commute in $G$.  But it is easily checked (using the fact that
all entries of $s(f_i) - f_i$ are divisible by $2$) that $s(f_1)^2$ and
$s(f_2)^2$ don't commute.
\vs

\textbf{Example 2}.  Let $m$ and $n$ be positive integers with $n \geqslant
3$, and let $p$ be a prime with $p^m > 3$.  Let $G$ (resp.\ $H$) be the
multiplicative group of invertible $n \times n$ matrices over the ring $\Z
/(p^{3m}\Z)$ (resp.\ $\Z /(p^{m}\Z)$), and let $\nu: G \to H$ be the
canonical surjection.  We write $I$ for the identity element in $G$ and in
$H$.  The kernel of $\nu$ is the group of matrices of the form $I + p^mf$
where $f$ is in $G$.  For any $i,j$ with $1 \leqslant i < j \leqslant n$ let
$\delta_{ij}$ be the $n \times n$ matrix which has $1$ in the $ij$-th place
and $0$ elsewhere; then $I + \delta_{ij}$ is an element of $G$ and $\nu(I +
\delta_{ij})$ has order $p^m$.  The liftings of $\nu(I + \delta_{ij})$ to
$G$ are the matrices of the form $I + \delta_{ij} + p^mf$ with $f$ in $G$.
Now we repeat a calculation from Evans, Hodges and Hodkinson \cite{evhoho:1}
Prop.\ 3.7.  The element $(I + \delta_{ij} + p^mf)^{p^m}$ is
\[ I + \binom{p^m}{1}(\delta_{ij} + p^mf) + \binom{p^m}{2}(\delta_{ij} +
p^mf)^2 + \binom{p^m}{3}(\delta_{ij} + p^mf)^3 + \ldots  \]
Since $\delta_{ij}\delta_{ij} = 0$, $p^{3m}x = 0$ in $\Z /(p^{3m}\Z )$ and
$p^m > 3$, this multiplies out to
\[  I + p^m\delta_{ij} + p^{2m}f + \frac{p^{2m}(p^{m}-1)}{2}(\delta_{ij}f
+ f\delta_{ij}) +
\frac{p^{2m}(p^{m}-1)(p^{m}-2)}{6}\delta_{ij}f\delta_{ij}.  \]
To apply these calculations to a concrete example, take
\[  g_1 = I + \delta_{12}, \ \ g_2 = I + \delta_{23} \]
in $H$, and let $s$ be a weak splitting of $\nu$.  Then
\[  s(g_1) = I + \delta_{12} + p^mf_1, \ \ s(g_2) = I + \delta_{23} +
p^mf_2  \]
for some $f_1$, $f_2$ in $G$.  Since $s$ is a weak splitting,
\[ s(g_1)^{p^m}s(g_2)^{p^m} = s(g_2)^{p^m}s(g_1)^{p^m}.\] 
But our calculations show at once that
\[  s(g_1)^{p^m}s(g_2)^{p^m} - s(g_2)^{p^m}s(g_1)^{p^m} =
p^{2m}\delta_{13} \neq 0. \]
This contradiction proves that $\nu$ doesn't weakly split.
\vs

Now suppose $\cons$ is a construction in the model $M$, and $B$
is a structure in the graph of \textbf{C}.  Let $A$ be $B^-$.  We
write $\Aut(A)$ and $\Aut(B)$ for the automorphism groups of $A$ and
$B$ respectively.  By (4) in the definition of constructions, each
automorphism $g$ of $B$ restricts to an automorphism $\nu_B(g)$ of $A$.
This map $\nu_B: \Aut(B) \to \Aut(A)$ is clearly a homomorphism;  by
(5) in the definition of constructions it is surjective.

We say that $B$ is \emph{(weakly) natural over} $A$ if the map $\nu_G$
(weakly) splits.  We say that the construction $\cons$ is \emph{(weakly)
natural} if for every $B$ in the graph of $\cons$, $\nu_B$ (weakly) splits.
(Our paper \cite{hosh:1} explained how this terminology connects with the
notion of a natural transformation.  In a related context Harvey Friedman
\cite{fr:1} used the term `naturalness' in a weaker sense.)
\vs

\textbf{Example 3}.  Let $G$ and $H$ be as in Example 1.  Since $n \times n$
upper triangular matrix groups are nilpotent of class $n-1$, $G$ is a finite
soluble group. So by Shafarevich \cite{sh:1} there is a Galois extension $K$
of the field $\mathbb Q$ of rationals such that $G$ is the Galois group of
$K/{\mathbb Q}$.  Let $k$ be the fixed field of the kernel $G_0$ of $\nu: G
\to H$.  Then $H$ is the Galois group of the extension $k/{\mathbb Q}$.  One
can write a set-theoretic description of these fields---up to
isomorphism---as a construction \textbf{C} where $K$ is in the graph and $k$
is picked out within $K$ by the relation symbol $P$.  This construction
\textbf{C} is small (in fact unitype) and $\emptyset$-representable in any
model of set theory, and it is not weakly natural.
\vs

\textbf{Example 4}.  Let $G$ and $H$ be as in Example 2.  Let $B$ (resp.\
$A$) be the direct sum of $n$ copies of the abelian group $\Z /(p^{3m}\Z
)$ (resp.\ $\Z /(p^m\Z )$), and identify $A$ with $p^{2m}B$.  Let the
relation symbol $P$ pick out $A$ within $B$.  Then $G$ (resp.\ $H$) is the
automorphism group of $B$ (resp.\ $A$), and $\nu: G \to H$ is the map
induced by restriction.  By the result of Example 2, the construction of
$B$ over $A$, which is again unitype and $\emptyset$-representable in any
model of set theory, is not weakly natural.
\vs

In \cite{hosh:1} we conjectured that there are models of set theory
in which each representable construction is uniformisable if and only
if it is natural.  See section \ref{se:7} below for some of the
background to this.  Section \ref{se:3} will show that no reasonable
version of this conjecture is true.  Sections \ref{se:4}--\ref{se:6}
will show that there are models in which uniformisability implies
weak naturality.  Section \ref{se:7} solves some of the problems
raised in \cite{ho:1} and \cite{hosh:1}.

\section{Uniformisability}   \label{se:3}

A structure $B$ is said to be \emph{rigid} if it has no nontrivial
automorphisms.  We shall say that a construction $\cons$ is
\emph{rigid-based} if for every structure $B$ in the graph of $\cons$, $B^-$
has no nontrivial automorphisms.  A rigid-based construction is trivially
natural.

Let $M$ be a transitive model of set theory.  We shall use a device that
takes any construction $\cons$ in $M$ to a construction $\cons^r$, called
its \emph{rigidification}.  Each structure $B^-$ in the domain is replaced
by a two-part structure $B^{r-}$, where the first part is $B^-$ and the
second part consists of the transitive closure of the set $P^B$ with a
membership relation $\varepsilon$ copying that in $M$.  Now $B^r$ is defined
to be the amalgam of $B$ and $B^{r-}$, so that $B^{r-}$ is $(B^r)^-$.  Then
$\cons^r$ is the closure of the class
\[ \{ B^r : B \textrm{ in the graph of } \cons \}  \]
under isomorphism in $M$.  It is clear that $\cons^r$ and the map
$B \mapsto B^r$ are definable in $M$ using no parameters beyond
those in the formulas representing $\cons$.  

\begin{lemma}  \label{le:1}
If $\cons$ is any construction, then $\cons^r$ is rigid-based, natural and
not small. 
\end{lemma}

\textbf{Proof}.  If $B^-$ is in the domain of $\cons$, then $B^{r-}$ is
rigid because its set of elements is transitively closed; so $\cons^r$ is
rigid-based.  Naturality follows at once.  Since the domain of $\cons$ is
closed under isomorphism, the relevant transitive closures are arbitrarily
large.  $\Box$

\begin{theorem} \label{th:2}
There is no transitive model of ZFC in which both the following
are true:
\begin{enumerate}
\item[(a)]  Every rigid-based construction in $M$ is uniformisable.
\item[(b)] Every unitype uniformisable construction in $M$ is weakly natural.
\end{enumerate}
In particular there is no transitive model of ZFC in which the
natural constructions are exactly the uniformisable ones.
\end{theorem}

\textbf{Proof}.  Suppose $M$ is a counterexample.  Let $\cons$ be
a unitype non-weakly-natural construction in $M$, such as Example
2 in section \ref{se:2}.  Then $\cons^r$ is rigid-based
and hence uniformisable by assumption.  

But we can use the uniformising formula of 
$\cons^r$ to uniformise $\cons$ with the same parameters.  So
by the assumption on $M$ again, $\cons$ is weakly natural;  
contradiction.  $\Box$
\vs

The next result gives some finer information about small constructions.

\begin{theorem}  \label{th:3}
Let $M$ be a transitive model of ZFC, $Y$ a set in $M$ and $\bar{c}$ a
well-ordering of $Y$ in $M$.  Assume:
\begin{quote}
In $M$, if $X$ is any set, then every unitype $X$-representable 
rigid-based construction is $X \cup Y$-uniformisable.  
\end{quote}
Then 
\begin{quote}
In $M$, every small $\emptyset$-representable construction is
$\{\bar{c}\}$-uniformisable,
\end{quote} 
and hence there are unitype $\{\bar{c}\}$-uniformisable constructions that
are not weakly natural.
\end{theorem}

\textbf{Proof}.  Let $\gamma$ be the length of $\bar{c}$.  Write $\bar{v}$
for the sequence of variables $(v_i : i < \gamma)$.  In $M$ we can
well-order (definably, with no parameters) the class of pairs $\langle
j,\psi \rangle$ where $j$ is an ordinal and $\psi(x,y,z,\bar{v})$ is a
formula of set theory.  We write $H_j$ for the set of sets hereditarily of
cardinality less than $\aleph_j$ in $M$.

Let $\cons$ be a small $\emptyset$-representable construction in $M$.  Then
$\cons^r$ is an $\emptyset$-representable rigid-based construction.  It is
not small; but if $B$ is any structure in the graph of $\cons$, let
$\cons_B$ be the construction got from $\cons^r$ by restricting the graph to
structures isomorphic to $B^r$.  Then $\cons_B$ is a unitype and
$\{B\}$-representable rigid-based construction, so by assumption it is
$\{B\} \cup Y$-uniformisable, say by a formula $\psi_B(-,-,B,\bar{c})$ where
$B, \bar{c}$ are the parameters.

By the reflection principle in $M$ there is an ordinal $j$ such that
\begin{quote}
$M \models \exists C (C \in \cons_B \land C^- = B^{r-} \land
C$ is the unique set such that ``$H_j \models \psi_B(B^{r-},C,B,
\bar{c})$'').
\end{quote}  
Hence in $M$ there is a first pair $\langle j_B, \psi_B \rangle$, 
definable from $B$, such that
\begin{quote}
$M \models \exists C (C \in \cons_B \land C^- = B^{r-} \land C$ is
the unique set such that ``$H_{j_B} \models \psi_B(B^{r-},C,B,
\bar{c})$'').
\end{quote}
Since all of this is uniform in $B$, it follows that the construction
$\cons$ is $\{\bar{c}\}$-uniformisable in $M$ by the formula
\[  y = C|L \textrm{ where } H_{j_B} \models \psi_B(B^{r-},C,B,
\bar{c}).  \]
The last clause of the theorem follows by choosing $\cons$ suitably,
for example as in Example 2 of section \ref{se:2}.  $\Box$

\section{The model} \label{se:4}

\begin{theorem}  \label{th:4}
There is a transitive model $N$ of ZFC in which every
uniformisable construction is weakly natural.
\end{theorem}

The next three sections are devoted to proving this theorem.  We use
forcing.  The central idea is to consider a construction $\cons$ whose
parameters lie in the ground model, and introduce a highly homogeneous
generic structure $B^\star$ into the graph of $\cons$; by homogeneity
$B^\star$ must be highly symmetrical over $B^{\star -}$.  Since the
parameters of a construction may lie anywhere in the set-theoretic universe,
we have to iterate this idea right up through the universe.  So we need to
build $N$ by a proper class iteration.

Our forcing notation is mainly as in Jech~\cite{je:1}.  Thus $p < q$ means
that $p$ carries more information than $q$.  We write $\dot{x}$ for a
boolean name of the element $x$ of $N$, and $\check{x}$ for the canonical
name of an element $x$ of the ground model.  If $y$ is a boolean name and
$G$ a generic set, we write $y[G]$ for the element named by $y$ in the
generic extension by $G$.  Our notion of forcing is of the kind described in
Menas~\cite{me:2} as `backward Easton forcing', and we shall borrow some
technical results from Menas' paper.

We start from a countable transitive model $M$ of ZFC + GCH. In $M$,
$\Lambda$ is a definable continuous monotone increasing function from
ordinals to infinite cardinals, with the property that for any ordinal
$\alpha$, $\Lambda(\alpha + 1) > \Lambda(\alpha)^{+}$.  Our notion of
forcing $\R_\infty$ will be defined by induction on the ordinals.  We start
with a trivial ordering $\R_0$.  At limit ordinals we take inverse limits.

For each ordinal $\alpha$ we shall define an $\R_\alpha$-name
$\dot{\nofS}_\alpha$; then $\R_{\alpha + 1}$ will be $\R_\alpha \otimes
\dot{\nofS}_\alpha$.  To define this name, let $\lambda$ be an infinite
cardinal.  We consider a notion of forcing, $\nof_\lambda$.  In
$\nof_\lambda$, conditions are partial maps $p: \lambda^+ \times \lambda^{+}
\times \lambda^+ \to 2$ with domain of cardinality $\leq \lambda$.  Write
$TP(\lambda)$ for a set-theoretical term which defines the notion of forcing
$\nof_\lambda$.  For each ordinal $\alpha$, we choose $\dot{\nofS}_\alpha$
to be an $\R_\alpha$-name such that $|| \dot{\nofS}_\alpha =
TP(\Lambda(\check{\alpha} + 1))||_{\R_\alpha} = 1$.

This defines a proper class notion of forcing, $\R_\infty =$ 
direct limit of $\langle\R_\alpha:\alpha\mbox{ ordinal}\rangle$.

\begin{lemma} \label{le:5}
For each ordinal $\alpha$, suppose $\Lambda(\alpha)$ is a cardinal
in $M^{\R_\alpha}$.  Then with $\R_\alpha$-boolean value 1,
$\dot{\nofS}_\alpha$ is $\Lambda(\alpha + 1)$-closed and satisfies
the $\Lambda(\alpha + 1)^{+}$-chain condition.
\end{lemma}

$\proof$.  Straightforward.  $\Box$

\begin{lemma}  \label{le:6}
\begin{enumerate}
\item[(a)] For every successor ordinal $\alpha$,
$\Lambda(\alpha)$ and $\Lambda(\alpha)^{+}$
are cardinals with $\R_\infty$-boolean value 1.
\item[(b)] For each successor ordinal $\alpha$, $|\R_{\alpha}| =
\Lambda (\alpha)^{+}$;  for each limit ordinal $\delta$,
$|\R_\delta| \leq \Lambda (\delta)^{+}$.
\item[(c)] For each ordinal $\alpha$, $\R_{\alpha}$ satisfies
the $\Lambda (\alpha)^{+}$-chain condition.
\end{enumerate}
\end{lemma}

\proof.  We prove all parts simultaneously by induction.  Suppose $R_\alpha$
has cardinality $\leq \Lambda (\alpha)^{+}$ and satisfies the
$\Lambda(\alpha)^{+}$-chain condition.  Put $\lambda = \Lambda (\alpha + 1)
> \Lambda(\alpha)^{+}$.  All cardinals $\geq \lambda$ are cardinals with
$R_\alpha$-boolean value 1.  Let $\dot{q}$ be an element of
$\dot{\nofS}_\alpha$.  Then $\dot{q}$ has cardinality $\leq \lambda$ with
$R_\alpha$-boolean value 1, and $\R_\alpha$ satisfies the $\lambda^+$-chain
condition, so with boolean value 1 the domain of $\dot{q}$ lies within some
$\gamma < \lambda^+$.  Now $\dot{p}$ can be taken to be a map from the set
$\gamma$ to the regular open algebra $RO(R_\alpha)$, which has cardinality
$\leq (|\R_\alpha|^{+})^{|\R_\alpha|} = |\R_\alpha|$.  The number of such
maps is at most $(|\R_\alpha|^{+})^\lambda = \lambda^+$.  Therefore
$|\R_\alpha \otimes \dot{\nofS}_\alpha| = \lambda^+$.  Also $R_{\alpha + 1}$
satisfies the $\lambda^{+}$-chain condition by Lemma~\ref{le:5}, using a
standard argument on iterated forcing (Menas~\cite{me:2} Proposition 10(i)).

Now with boolean value 1, $\dot{\nofS}_\alpha$ is $\lambda$-closed and
satisfies the $\lambda^{+}$-chain condition, so cardinals are preserved in
passing from $\R_\alpha$ to $\R_{\alpha + 1}$.

We turn to limit ordinals $\delta$.  The cardinality of $\R_\delta$ is at
most the product of the cardinalities of $\R_\alpha$ with $\alpha < \delta$,
hence at most $\Lambda(\delta)^+$.  It follows at once that $\R_\delta$
satisfies the $\Lambda(\delta)^{+}$-chain condition and preserves all
cardinals from $\Lambda(\delta)^{+}$ upwards.

It remains to show that for successor ordinals $\alpha + 1$, the cardinals
$\Lambda(\alpha + 1)$ and $\Lambda(\alpha + 1)^{+}$ are not collapsed by
$\R_\infty$.  Using the next lemma (which doesn't depend on the clause we
are now proving), $\R_\infty$ can be written as $\R_\alpha \otimes
\dot{\nofS}_{\alpha + 1} \otimes \dot{\R}_{\alpha + 1,\infty}$.  The first
factor satisfies the $\Lambda(\alpha + 1)$-chain condition and hence
preserves these two cardinals.  The third factor preserves them since it is
$\Lambda(\alpha + 1)^{+}$-closed with boolean value 1.  The middle factor is
$\Lambda(\alpha + 1)$-closed with boolean value 1, so that it preserves
$\Lambda(\alpha + 1)$ and $\Lambda(\alpha + 1)^{+}$.  $\Box$

\begin{lemma} \label{le:7}
For each ordinal $\alpha$ there is a proper class notion of
forcing $\dot{\R}_{\alpha,\infty}$ such that
\begin{enumerate}
\item[(a)] $\R_\infty$ is isomorphic to 
$\R_\alpha \otimes \dot{\R}_{\alpha,\infty}$;
\item[(b)] In $M^{\R_\alpha}$, $\dot{\R}_{\alpha,\infty}$ is the direct 
limit of iterated notions of forcing $\dot{\R}_{\alpha,\beta}$
in such a way that for each $\beta > \alpha$, $\R_\beta$ is
isomorphic to $\R_\alpha \otimes \dot{\R}_{\alpha,\beta}$;
\item[(c)] For each successor ordinal $\alpha$, 
$\dot{\R}_{\alpha,\infty}$ 
is $\Lambda(\alpha + 1)$-closed with $M^{\R_\alpha}$-boolean
value 1.
\end{enumerate}
\end{lemma}

\proof.  As Menas~\cite{me:2} Propositions 11 and 10(i),
using the previous lemma.  $\Box$
\vs

The model $N$ for the theorem will be an $R_\infty$-generic
extension of $M$.  

\begin{lemma}  \label{le:8}
$N$ is a model of ZFC.
\end{lemma}

\proof.  Menas~\cite{me:2} Proposition 14 derives this from the
previous lemma. $\Box$

\begin{lemma}  \label{le:9}
If $\dot{x}$ is an $\R_\infty$-name of a subset of $\Lambda(\alpha + 1)$,
and $r \in \R_\infty$, then there are $s \in \R_\infty$ and an
$\R_\alpha$-name $\dot{y}$ such that $s \leq r$ and $s \forces_{\R_\infty}$
``$\dot{x} = \dot{y}$''.
\end{lemma}

\proof.  This follows from the fact that $\dot{R}_{\alpha, \infty}$
is $\Lambda(\alpha + 1)$-closed.  $\Box$

\begin{lemma}  \label{le:10}
Let $\alpha$ be an ordinal.  Then if $\lambda$ is $\Lambda(\alpha + 1)$ or
$\Lambda(\alpha + 1)^{+}$, we have $2^\lambda = \lambda^+$ in $N$.
\end{lemma}

$\proof$.  Suppose $\lambda = \Lambda(\alpha + 1)$.  Then $2^\lambda \leq
\mu$ with $R_\alpha$-boolean value 1, where $$\mu = |RO(\R_\alpha)|^\lambda
\leq (|\R_\alpha|^{\Lambda(\alpha)^{+}})^\lambda = \lambda^+.$$ With
$\R_\alpha$-boolean value 1, $\R_{\alpha,\infty}$ is $\lambda$-closed and
hence adds no new subsets of $\lambda$.  Similar calculations apply to the
other cases.  $\Box$
\vs

A notion of forcing $\R$ is said to be {\em homogeneous} if for any two
elements $p, q$ of $\R$ there is an automorphism $\sigma$ of $\R$ such that
$\sigma (p)$ and $q$ are compatible.

\begin{lemma} \label{le:11}
The notion of forcing $\R_\infty$ is homogeneous.
\end{lemma}

\proof.  Menas~\cite{me:2} Proposition 13 proves this under the assumption
that each step of the iteration is homogeneous with boolean value 1.  That
assumption holds here.  $\Box$
\vs

If $\alpha$ is an automorphism of the notion of forcing $\R$, then $\alpha$
induces an automorphism $\alpha^\star$ of the boolean universe $M^P$.  Also
$\alpha$ takes any $\R$-generic set $G$ over $M$ to the $\R$-generic set
$\alpha G$.

\begin{lemma}  \label{le:12}
For every element $\dot{x}$ of $M^\R$ we have
\[\dot{x}[G] = (\alpha^\star)\dot{x}[\alpha G]\]
\end{lemma}

\proof.  Immediate.  $\Box$
\vs

To save notation we write $\alpha^\star$ as $\alpha$.  We note that
$(\alpha\beta)^\star = \alpha^\star\beta^\star$, which removes one possible
source of ambiguity.

\section{The generic copies of $A$, $B$} \label{se:5}
\vs

As explained earlier, our model $N$ in the theorem will be $M[G]$ where $G$
is an $\R_\infty$-generic class over $M$.  Henceforth $\cons$ is a
construction which is uniformisable in $N$ with uniformising formula $\phi$;
we want to show that $\cons$ is weakly natural.  Let $B$ be any structure in
the graph of $\cons$.  At the cost of adding $B$ as a parameter, we can
assume without loss that $\cons$ is unitype and its graph consists of
structures isomorphic to $B$.  We put $A = B^{-}$.  Choose an ordinal
$\alpha$ so that $A, B$ and the parameters of the formulas representing
$\cons$ all lie in $M[G \cap \R_\alpha]$, and $B, \mbox{Aut}(B)$ both have
cardinality $\leq \Lambda(\alpha + 1)$.  We can decompose $N$ as a two-stage
extension $M[G_\alpha][G_{\alpha,\infty}]$, where $G_\alpha = G \cap
\R_\alpha$ and $G_{\alpha,\infty}$ is $\dot{\R}_{\alpha,\infty}[
G_\alpha]$-generic over $M[G_\alpha]$.

At this point we adjust our notation.  We put $\lambda = \Lambda(\alpha +
1)$, and we rename $M[G_\alpha]$ as $M$.  By Lemma \ref{le:6}, $\lambda$ and
$\lambda^+$ in the old $M$ are still cardinals in the new $M$.  By
Lemma~\ref{le:6}, $N$ is constructed from the new $M$ by an iterated forcing
notion $\R_{\alpha,\infty} = \dot{\R}_{\alpha,\infty}[G_\alpha]$ with the
same properties as the forcing notion $\R_\infty$, with two differences.
First, the function $\Lambda$ must now be replaced by the function
$\Lambda_\alpha$ where $\Lambda_\alpha(\beta) = \Lambda(\alpha + \beta)$.
Second, $M$ need not satisfy the GCH everywhere; but this never matters.
(In fact it would be possible to make the GCH hold in the new $M$ and in
$N$, by adding extra factors in $\R_\infty$ to collapse the cardinalities of
power sets.)  One can check that all the preliminary lemmas \ref{le:5} to
\ref{le:11} still hold for this notion of forcing $\R_{\alpha,\infty}$.

We now write $\nof, \dot{\Q}$ for 
$\dot{\nofS}_\alpha[G_\alpha],
\dot{\R}_{\alpha + 1,\infty}[G_\alpha]$ respectively.  Thus
\[\R_{\alpha,\infty} = \gnof. \]  
We shall not need to refer to $G_\alpha$ again, and so we start afresh with
our notation for generic sets.

We shall write $N$ as $M[G_1][G_2]$ where $G_1$ is $\nof$-generic over $M$
and $G_2$ is $\dot{\Q}[G_1]$-generic over $M[G_1]$.

We shall write $G$ for the $\gnof$-generic set $G_1 \otimes \dot{G}_2$ over
$M$, so that $N = M[G]$.  If $x$ is an element of $N$, we write $\dot{x}$
for a boolean name for $x$ in the forcing language for $\gnof$.  Note that
every $\nof$-name over $M$ can be read as a $\gnof$-name too, so that there
is no need for a separate symbol for $\nof$-names.

The set $\bigcup G_1$ is a total map from $\lambda^+ \times \lambda^{+}
\times \lambda^+$ to 2.  For each $i < \lambda^+$ and $j < \lambda^{+}$, we
define $a_{ij} = \{ k < \lambda^+ : \bigcup G_1 (i,j,k) = 1 \}$ and $a'_i =
\{a_{ij} : j < \lambda^{+}\}$, so that $a'_i$ is a set of $\lambda^{+}$
independently generic subsets of $\lambda^+$.  If $a$ and $b$ are (in
$M_1[G_2]$) sets of subsets of $\lambda^+$, we put $a \equiv b$ iff the
symmetric difference of $a$ and $b$ has cardinality $\leq \lambda$.  We
write $a_i$ for the equivalence class $(a'_i)^\equiv$.  The boolean names
$\dot{a}_{ij}, \dot{a}_i$ can be chosen in $M^{\nof}$ independently of the
choice of $G$.

Consider again the structures $A$ and $B$ in $M$.  Without loss we can
suppose that dom$(A)$ is an initial segment of $\lambda$.  In $M[G_1]$ there
is a map $e$ which takes each element $i$ of $A$ to the corresponding set
$a_i = \dot{a}_i[G_1]$; by means of $e$ we can define a copy $A^*$ of $A$
whose elements are the sets $a_i$ ($i \in \mbox{dom}(A)$).  Again the
boolean names $\dot{A}^*, \dot{e}$ can be chosen to be independent of the
choice of $G$.

Since $A, B$ and the parameters of the uniformising formula $\phi$ lie in
$M$, and the notion of forcing $\gnof$ is homogeneous by Lemma~\ref{le:11},
the statement ``$\phi$ defines a construction on the class of structures
isomorphic to $A$, which takes $A$ to $B$'' is true in $N$ independently of
the choice of $G$.  In particular there are $\gnof$-boolean names
$\dot{B}^*, \bnep$ such that
\begin{eqnarray}  \label{eq:2}
& & || \dot{B} \mbox{ is the unique structure such that
} \phi(\dot{A}^*, \dot{B}^*) \mbox{ holds},\\
& & \dot{e} : \check{A} \to \dot{A}^* \mbox{ is
the isomorphism such that }
\dot{e}(\check{\imath}) = \dot{a}_i \mbox{ for}
\nonumber \\
& & \mbox{each } i \in \mbox{dom}(\check{A}), 
\mbox{ and }
\bnep : \check{B} \to \dot{B}^* \mbox{ is 
an isomorphism}
\nonumber \\
& & \mbox{which extends }
\dot{e}||_{\gnof} = 1.  \nonumber
\end{eqnarray}

\begin{lemma} \label{le:13} 
Let $G$ be $\gnof$-generic over $M_1$.  Then:

(a)  {\rm Aut}$(A)^M$ =  
{\rm Aut}$(A)^{M[G_1]}$ = {\rm Aut}$(A)^{M[G]}$.

(b)  {\rm Aut}$(B)^M$ =  
{\rm Aut}$(B)^{M[G_1]}$ = {\rm Aut}$(B)^{M[G]}$.

(c)  The set of maps from {\rm Aut}$(A)$ to {\rm Aut}$(B)$ is the same
in $M$ as it is in $M[G_1]$ and $M[G]$.
\end{lemma}

\proof.  Using Lemma~\ref{le:5} and Lemma~\ref{le:6}(c), $\nof$ is
$\lambda$-closed over $M$, and $\dot{\Q}$ is $\lambda$-closed over $M[G_1]$.
Hence no new permutations of $A$ or $B$ are added since $|A| \leq |B| \leq
\lambda$ in $M$; this proves (a), (b).  Likewise (c) holds since
$|\mbox{Aut}(A)| \leq |\mbox{Aut}(B)| \leq \lambda$ in $M$.  $\Box$
\vs

We regard $\Aut (A)$ as a permutation group on $\lambda^+$ by letting it fix
all the elements of $\lambda^+$ which are not in dom$(A)$.

By a {\em neat map} we mean a map $\alpha : \lambda^{+} \to \Aut (A)$ in $M$
which is constant on a final segment of $\lambda^{+}$; we write $\neat$ for
the set of neat maps.  We write $\pi$ for the map from $\neat$ to $\Aut (A)$
which takes each neat map to its eventual value.  We write $\neat^-$ for the
set of all neat maps $\alpha$ with $\pi(\alpha) = 1$.  For each ordinal $i <
\lambda^{+}$ we write $\neat_i$ for the set of neat maps $\alpha$ such that
$\alpha(j) = 1$ for all $j < i$.  We write $\neat^-_i$ for $\neat^- \cap
\neat_i$.

We can regard $\alpha$ as a
permutation of the set $\lambda^+ \times \lambda^{+}
\times \lambda^+$ by putting
\[\alpha(i,j,k) = (\alpha(j)(i),j,k).\]
Then $\alpha$ induces an automorphism of $\nof$.

\begin{lemma}  \label{le:14}
If $\alpha$ and $\beta$ are distinct neat maps then they induce
distinct automorphisms of $\nof$.  Identifying each neat map with
the automorphism it induces, $\neat$ forms a group with 
subgroups $\neat^-$, $\neat_i$ $(i < \lambda^+)$;
the map $\pi : \neat \to \Aut (A)$ is
a group homomorphism.
\end{lemma}

\proof.  From the definitions.  $\Box$
\vs

The automorphism $\alpha$ can be extended to an automorphism of $\gnof$ in
many different ways, by induction on $\dot{\R}$ as an iterated notion of
forcing.  Each factor of $\dot{\R}$ is with boolean value 1 the set of all
maps from $X$ to 2 of cardinality $\leq \mu$, where $X$ is $\mu^+ \times
\mu^{+} \times \mu^+$ for some cardinal $\mu$.  If $\Sigma '$ is (in $M$)
the group of permutations of $X$, then an automorphism of the factor of
$\dot{\R}$ is determined by an element of $\Sigma '$ and a permutation of
the boolean values.  For each ordinal $i$ let $\Sigma_i$ be in $M_1$ the
product of the permutation groups $\Sigma '$ for the first $i$ factors of
$\dot{\Q}$, and let $\Sigma$ be the direct limit of the $\Sigma_i$ in $M_i$.
Then an automorphism $\alpha$ of $\nof$ and an element $\sigma$ of $\Sigma$
together determine an automorphism $\langle \alpha , \sigma \rangle$ of
$\gnof$, and hence of $M^{\gnof}$.

\begin{lemma}  \label{le:15}
The actions of the group $\cal N$ of neat maps and the group $\Sigma$ on
$\gnof$ commute with each other.
\end{lemma}

\proof.  The class $\gnof$ is $\nof \times \hat{\Q}$ where $\hat{\Q}$ is a
class of boolean-valued subsets of a class $X$ which is definable in $M$;
the action of $\Sigma$ is through its action on $X$.  Thus each element of
$\hat{\Q}$ is essentially a set of ordered pairs $\langle x , y \rangle$
where $x \in X$ and $y \in \nof$. Since $\Sigma$ and $\cal N$ act
respectively on the first and second coordinates, the actions on $\hat{\Q}$
commute. The group $\Sigma$ keeps $\nof$ fixed.  $\Box$

\begin{lemma} \label{le:16} 
Suppose $\alpha : \lambda^{+} \to \Aut (A)$ is neat and $\alpha '$ is an
automorphism of $\gnof$ extending $\alpha$. Then the action of $\alpha '$ on
$M^{\gnof}$ setwise fixes the set $\{ \dot{a}_i : i \in {\rm dom}(A) \}$ of
canonical names of the elements of $\dot{A}^*[G]$, and it acts on this set
in the way induced by $\pi(\alpha)$ and the map $i \mapsto \dot{a}_i$.  Thus
$\alpha '(\dot{a}_i) = \dot{a}_{\pi(\alpha)(i)}$.
\end{lemma}

\proof.  Write out the names!  (They lie in $M^{\nof}$, so that the
extension from $M^\nof$ to $M^{\gnof}$ is irrelevant.)  $\Box$
\vs

If $G$ is $\gnof$-generic over $M$, then so is $\langle \alpha, \sigma
\rangle G$ for every neat map $\alpha$ and every $\sigma \in \Sigma$, since
$\alpha, \sigma \in M$.

\begin{lemma}  \label{le:17}
For each element $i$ of $A$, each neat map $\alpha$ and each $\sigma \in
\Sigma$, $\dot{a}_{\pi(\alpha)(i)}[\langle \alpha, \sigma \rangle G] =
\dot{a}_i[G]$.  In particular $\dot{A}^*[\langle \alpha, \sigma \rangle G] =
\dot{A}^*[G]$.
\end{lemma}

$\proof$.  By Lemma~\ref{le:16}, $\dot{a}_{\pi(\alpha)
(i)}[\langle \alpha, \sigma \rangle G] = 
(\alpha \dot{a}_i)[\langle \alpha, \sigma \rangle G]$.  Then by Lemma
\ref{le:12} and the fact that $\alpha \dot{a}_i$ lies in $M^{\nof}$,
\[  (\alpha \dot{a}_i)[\langle \alpha, \sigma \rangle G] = 
(\alpha \dot{a}_i)[\alpha  G_1] = \dot{a}_i[G_1] = \dot{a}_i[G]. 
\]

$\Box$
\vs

We write $\bnep^{-1}$ for a boolean name such that $\bnep^{-1}[G] =
(\bnep[G])^{-1}$ for all generic $G$.

\begin{lemma} \label{le:18}
Suppose $\alpha$ is a neat map, $\sigma \in \Sigma$ and $G$ is
$\gnof$-generic over $M_1$.  Then $\dot{B}^*[\langle \alpha, \sigma
\rangle^{-1} G]$ $= \dot{B}^*[G]$, and the map $(\bnep ^{-1} \circ \langle
\alpha, \sigma \rangle \bnep)[G]$ is an automorphism of $B$ which extends
$\pi(\alpha)$.
\end{lemma}

\proof.  Since $M[\langle \alpha, \sigma \rangle^{-1} G] = M_1[G]$ and
$\dot{A}^*[\langle \alpha, \sigma \rangle^{-1} G] = \dot{A}^*[G]$,
(\ref{eq:2}) (before Lemma \ref{le:13}) tells us that $\dot{e} [\langle
\alpha, \sigma \rangle^{-1} G] (i) = \dot{a}_i[\langle \alpha, \sigma
\rangle^{-1} G]$ for each $i \in \mbox{dom}(A)$, and that $\dot{B}^*[\langle
\alpha, \sigma \rangle^{-1} G] = \dot{B}^*[G]$ and $\bnep [G]^{-1}\circ
\bnep [\langle \alpha, \sigma \rangle^{-1} G]$ extends $\dot{e}[G]^{-1}\circ
\dot{e}[\langle \alpha, \sigma \rangle^{-1} G]$.  Now using Lemma
\ref{le:17},
\[ \dot{e}[G]^{-1}\circ\dot{e}[\langle \alpha, \sigma \rangle^{-1} G](i) 
= \dot{e}[G]^{-1}(\dot{a}_i[\langle \alpha, \sigma \rangle^{-1} G])\] 
\[= \dot{e}[G]^{-1}(\dot{a}_{\pi(\alpha)(i)}[G]) = \pi(\alpha)(i).\]
$\Box$
\vs

\begin{lemma} \label{le:19}
For every neat map $\alpha$, each $\sigma \in \Sigma$ and all $\langle
p,\dot{q} \rangle \in \gnof$ there are $\langle p',\dot{q}' \rangle
\leqslant \langle p,\dot{q} \rangle $ and $g \in \Aut{B}$ such that
\[  \langle p',\dot{q}' \rangle  \forces_{\gnof} 
\sigma (\bnep^{-1}) \circ \alpha\sigma (\bnep) = \check{g}.  \]
\end{lemma}

\proof.  Let $f$ be $\pi(\alpha)$.
By Lemma \ref{le:18} we have
\[|| \sigma\bnep^{-1} \circ \alpha\sigma \bnep \textrm{ is an automorphism 
of } B \textrm{ extending } \check{f}||_{\gnof} = 1.\]
Unpacking the existential quantifier in ``an automorphism of $B$''
gives the lemma.  $\Box$
\vs

Consider any $i < \lambda^{+}$.  Given $\langle p,\dot{q} \rangle \in
\gnof$, define $t_{p,\dot{q},i}$ to be the set of all triples
$(f,g,\sigma)$, with $f \in \Aut (A)$, $g \in \Aut (B)$ and $\sigma \in
\Sigma$, such that for some $\alpha \in \neat_i$, $\pi(\alpha) = f$ and
\[\langle p,\dot{q} \rangle \forces_{\gnof}\sigma(\bnep^{-1}) \circ
\alpha\sigma(\bnep) = \check{g}.\] Clearly if $\langle p',\dot{q}' \rangle
\leqslant \langle p,\dot{q} \rangle$ then $t_{p',\dot{q}',i} \supseteq
t_{p,\dot{q},i}$.  Since there are only a set of values for $\sigma(\bnep)$
and $\sigma(\bnep^{-1})$ with $\sigma \in \Sigma$.  it follows that there is
$\langle p_i,\dot{q}_i \rangle$ such that for all $\langle p',\dot{q}'
\rangle \leqslant \langle p,\dot{q} \rangle$,
\[t_{p',\dot{q}',i} = t_{p_i,\dot{q}_i,i}.\]
We fix a choice of $p_i, \dot{q}_i$, and we write $t_i$ for the resulting
value $t_{p_i,\dot{q}_i,i}$.  If $(f,g,\sigma)$ is in $t_i$, we write
$\alpha^i_{f,g,\sigma}$ for some $\alpha \in \neat_i$ such that
\[\langle p_i,\dot{q}_i \rangle \forces_{\gnof} \sigma(\bnep^{-1}) \circ 
\alpha\sigma(\bnep) = \check{g} \]
and  $\pi(\alpha) = f$.

\begin{lemma} \label{le:20}
For each $i < \lambda^{+}$, $t_i$ is a subclass 
of $\Aut (A) \times \Aut (B) \times \Sigma$ such that
\begin{enumerate}
\item[(a)] for each $(f,g,\sigma)$ in $t_i$, $g|A = f$;
\item[(b)] for each $f$ in $\Aut (A)$ and $\sigma$ in $\Sigma$
there is $g$ with $(f,g,\sigma)$ in $t_i$.
\end{enumerate}
(So $t_i(-,-,\sigma)$ is a first attempt at a lifting of the restriction map
from $\Aut (B)$ to $\Aut (A)$.)
\end{lemma}

\proof.  By Lemma \ref{le:19}.  $\Box$
\vs

We write $t^-_{p,\dot{q},i}$ for the set of pairs $(g,\sigma)$ such that
$(1,g,\sigma)$ is in $t_{p,\dot{q},i}$.  We write $\alpha^i_{g,\sigma}$ for
$\alpha^i_{1,g,\sigma}$; note that $\alpha^i_{g,\sigma}$ is in $\neat^-_i$
by Lemma \ref{le:18}.

\begin{lemma} \label{le:21}
For each $i < \lambda^{+}$ there are $\sigma_i$ in $\Sigma$, a condition
$p'_i \leqslant p_i$ and a boolean name $\dot{v_i}$ such that
\begin{enumerate}
\item[(a)] for each $i < \lambda^{+}$, $p'_i \forces_{\nof}
\dom(\sigma_i^{-1}\dot{q}_i) \subseteq \dot{v_i}$;
\item[(b)] for all $i < j < \lambda^{+}$, $||\dot{v}_i \cap
\dot{v}_j = \emptyset||_{\nof} = 1$;
\item[(c)] for all $i < j < \lambda^{+}$,
$\sigma_i\sigma_j = \sigma_j\sigma_i$.
\end{enumerate}
\end{lemma}

\proof.  By induction on $i < \lambda^{+}$.  As we choose the $p'_i$,
$\sigma_i$ and $\dot{v_i}$, we also choose an eventually zero sequence of
ordinals $\gamma_{\mu,i} < \mu^+$ in $M_1$ so that
\[||\dot{v_i} \subseteq \prod_{\mu} (\gamma_{\mu,i} \times\mu^+
\times \mu^+)||_{\nof} = 1.\]
Then when we have made our choices 
for all $i < j$, we first extend $p_i$ to $p'_i$ forcing the domain of
$\dot{q_i}$ to lie within some set 
\[X = \prod_{\mu < \mu'} (\gamma'_{\mu} \times \gamma'_{\mu}
\times \mu^{+}) \] lying in $M_1$, and we choose $\dot{w_i}$ to be a
canonical boolean name for this set $X$.  Then we choose $\sigma_i$ so that
$\sigma_i^{-1}$ moves $X$ to
\[\prod_{\mu<\mu'}\left( \left[ \bigcup_{j<i}\gamma_{\mu,i},\bigcup_{j<i}
\gamma_{\mu,i}+\gamma'_{\mu} \right)\times \gamma'_{\mu} \times \mu^{+}
\right) , \] 
(the product of products of three intervals), and we put $\dot{v_i} =
\sigma_i^{-1}\dot{w_i}$ and $\gamma_{\mu,j}=\bigcup_{j<i}\gamma_{\mu,i} +
\gamma'_{\mu}$.  $\Box$
\vs

We fix the choice of $\sigma_i$ and $\dot{v_i}$ $(i < \lambda^+)$ given by
this lemma.  Without loss we extend the conditions $p_i$ to be equal to
$p'_i$.

\begin{lemma} \label{le:22}
There is a stationary subset $S$ of $\lambda^+$ such that:
\begin{enumerate}
\item[(a)] for each $i \in S$ and $j < i$, 
the domain of $p_i$ is a subset of $i \times \dom A$;
\item[(b)] for each $i \in S$ and $j < i$, every map 
$\alpha^i_{f,g} : \lambda^{+} \to \Aut (A)$ 
is constant on $[i,\lambda^{+})$;
\item[(c)] for all $i, j \in S$, 
\[\{(f,g):(f,g,\sigma_i)\in t_i\}=\{(f,g):(f,g,\sigma_j)\in t_j\};\]
\item[(d)] there is a condition $p^{\star} \in \nof$ such that for all
$i \in S$, $p_i|(i \times \dom A) = p^{\star}$.
\end{enumerate}
\end{lemma}

$\proof$.  First, there is a club $C \subseteq \lambda^{+}$ on which (a) and
(b) hold.  Then by F\H{o}dor's lemma there is a stationary subset $S$ of $C$
on which (c) and (d) hold.  $\Box$

\section{The weak lifting}  \label{se:6}

In this section we use the notation $S$, $\sigma_i$, $\dot{v}_i$,
$p^{\star}$ from Lemmas \ref{le:21} and \ref{le:22}.  We write $s$ for the
constant value of
\[ \{(f,g) : (f,g,\sigma_i) \in t_i\} \ (i \in S) \]
from clause (c) of Lemma \ref{le:22}, and $s^-$ for the set of all $g$ such
that $(1,g) \in s$.  We write $\nu : \Aut(B) \to \Aut(A)$ for the
restriction map.

\begin{lemma}  \label{le:23}
The relation $s$ is a subset of $\Aut(A) \times \Aut(B)$ that projects onto
$\Aut(A)$, and if $(f,g)$ is in $s$ then $\nu(g) = f$.
\end{lemma}

\proof.  Lemma \ref{le:20}.  $\Box$
\vs

\begin{lemma} \label{le:24}
If $(f_1,g_1)$ and $(f_2,g_2)$ are both in $s$ then $(f_1f_2,g_1g_2)$ is in
$s$.
\end{lemma}

$\proof$.  In this and later calculations we freely use the fact (Lemma
\ref{le:15}) that the actions of $\mathcal N$ and $\Sigma$ on $\gnof$
commute.  Take any $i, j \in S$ with $i < j$.

Put $\alpha_1 = \alpha^j_{f_1,g_1,\sigma_j}$, $\alpha_2 =
\alpha^i_{f_2,g_2,\sigma_i}$ and $\alpha_3 = \alpha_1\alpha_2$.  Note that
$\alpha_1\alpha_2$ is in $\neat_i$ since $i < j$.
\vs

We have
\[ \langle p_j,\dot{q}_j \rangle \forces 
\sigma_j\bnep^{-1} \circ \alpha_3 \sigma_j (\bnep) = 
\sigma_j\bnep^{-1} \circ \alpha_1 \sigma_j (\bnep) \circ
(\sigma_j \alpha_1(\bnep))^{-1} \circ \alpha_3 \sigma_j (\bnep)  \]
and by assumption
\[  \langle p_j,\dot{q}_j \rangle  \forces 
\sigma_j\bnep^{-1} \circ \alpha_1 \sigma_j(\bnep) = \check{g_1}.  \]
So
\[ \langle p_j,\dot{q}_j \rangle  \forces 
\sigma_j\bnep^{-1} \circ \alpha_3 \sigma_j (\bnep) =
\sigma_j\check{g_1} \circ \sigma^j(\alpha_1(\bnep))^{-1} \circ 
\alpha_1(\alpha_2 \sigma_j \bnep). \]
Also by assumption
\[  \langle p_i,\dot{q}_i \rangle  \forces 
\sigma_i\bnep^{-1} \circ \alpha_2 \sigma_i(\bnep) = \sigma_i\check{g_2}.\]
Acting on this by $\alpha_1\sigma_j\sigma_i^{-1}$ gives
\[ \langle \alpha_1p_i,\alpha_1\sigma_j\sigma_i^{-1}\dot{q}_i\rangle 
\forces \alpha_1 \sigma_j(\bnep^{-1}) \circ \alpha_1\alpha_2\sigma_j\bnep =
\alpha_1\sigma_j\check{g_2}. \] 
Now $g_2$ is in the ground model and hence $\alpha_2\sigma_j\check{g_2} =
\check{g_2}$.  Also $\alpha_1p_i = p_i$ since the support of $p_i$ lies
entirely below $j$, and $\alpha_1 = \alpha^j_{g_1,\sigma_j}$ is the identity
in this region since it lies in $\neat_j$.  So we have shown that
\[ \langle p_i, \alpha_1\sigma_j\sigma_i^{-1}\dot{q}_i \rangle 
\forces \alpha_1\sigma_j\bnep^{-1} \circ 
\alpha_1\alpha_2\sigma_j\bnep = \check{g_2}.  \]

Now we note that $p_i \cup p_j$ is a condition in $P$, by (a), (d) of Lemma
\ref{le:22}.  Also $p_i \cup p_j$ forces that $\dom
(\sigma_i^{-1}\dot{q}_i)$ is disjoint from $\dom (\sigma_j^{-1}\dot{q}_j)$
by Lemma \ref{le:21}, and hence also that $\dom
\sigma_j\sigma_i^{-1}\dot{q}_i$ is disjoint from $\dom r_j$. From the action
of neat maps on $\Q$, $\dom \alpha_1\sigma_j\sigma_1^{-1}= \dom
\sigma_j\sigma_i^{-1}$.  This shows that $\langle p_i,\sigma_j\sigma_i^{-1}
\dot{q}_i \rangle$ and $\langle p_j,\dot{q}_j \rangle$ have a common
extension $\langle p',\dot{q}' \rangle$.  Putting everything together, we
have that
\[\langle p',\dot{q}'\rangle\forces\sigma_j\bnep^{-1}\circ\alpha_3\sigma_j
\bnep = \check{g_1}\check{g_2}.\] 
Since $\alpha_3$ is in $\neat_i$, this shows that
\[  (f_1f_2,g_1g_2) \in t_{p',\dot{q}',j}.  \]
Then by the maximality property of $\langle p_j,\dot{q}_j \rangle$,
\[  (f_1f_2,g_1g_2,\sigma_j) \in t_{p_j,j}  \]
so that $(f_1f_2,g_1g_2)$ is in $s$.
$\Box$
\vs

\begin{lemma}  \label{le:25}
If $g_1$ and $g_2$ are in $s^-$ then $g_1g_2 = g_2g_1$.
\end{lemma}
\vs

\proof.  Apply the proof of Lemma \ref{le:24} to $(1,g_1)$ and $(1,g_2)$.
In the notation of that proof, $\alpha_1$ commutes with $\alpha_2$,
$\alpha_1$ is the identity below $j$ and $\alpha_2$ is the identity below
$j$ (since $i, j \in S$.  But also $g_2$ lies in $s^-$, and this tells us
that $\alpha_2$ is the identity on $[j,\lambda^+)$.  In particular
$\alpha_1$ commutes with $\alpha_2$.

We follow the proof of Lemma \ref{le:24} but with $g_1$ and $g_2$
transposed, starting from the observation that
\[ \langle p_i,\dot{q}_i \rangle \forces\sigma_i\bnep^{-1} \circ
\alpha_3\sigma_i\bnep = \sigma_i \bnep^{-1} \circ \alpha_2\sigma_i \bnep
\circ \alpha_2\sigma_i\bnep^{-1} \circ \alpha_3\sigma_i\bnep. \] 
As before, we have that
\[  \langle p_i,\dot{q}_i \rangle  \forces 
\sigma_i\bnep^{-1} \circ \alpha_2\sigma_i\bnep = \check{g_2}  \]
and
\[\langle \alpha_2p_j,\alpha_2\sigma_i\sigma_j^{-1}\dot{q}_j \rangle
\forces \alpha_2\sigma_i\bnep^{-1} \circ \alpha_2\alpha_1\sigma_i\bnep =
\alpha_3\sigma_i\check{g_1}, \] 
recalling that $\alpha_1$ commutes with $\alpha_2$.  Now the support of
$p_j$ lies below $i$ or within $[j,\lambda^+)\times \dom A$, and $\alpha_2$
is the identity in both these regions, and so $\alpha_2(p_j) = p_j$.  Also
$p_i \cup p_j$ forces that $\dot{q_i}$ and
$\alpha_2\sigma_i\sigma_j^{-1}\dot{q_j}$ have disjoint domains.  So as
before there is $\langle p',\dot{q}' \rangle \leqslant \langle p_i,\dot{q}_i
\rangle$ and $\leqslant \langle p_j, \alpha_2\sigma_i\sigma_j^{-1}q_j
\rangle$ such that
\[  \langle p',\dot{q}' \rangle \forces 
\sigma_i\bnep^{-1} \circ \alpha_3\sigma_i\bnep  =
\check{g_2} \check{g_1}.  \]
As before, it follows that
\[ \langle p_i,\dot{q}_i \rangle \forces 
\sigma_i\bnep^{-1} \circ \alpha_3\sigma_i\bnep  =
\check{g_2} \check{g_1},  \]
and so
\[ \langle p_i, \sigma_j\sigma_i^{-1}\dot{q}_i \rangle \forces
\sigma_j\bnep^{-1} \circ \alpha_3\sigma_j\bnep  =
\check{g_2} \check{g_1}.  \]

Again there is a condition $\langle p'', \dot{q}'' \rangle \leqslant
\langle p_i, \sigma_j\sigma_i^{-1}\dot{q}_i \rangle$ and 
$\leqslant \langle  p_j,\dot{q}_j \rangle$.  Recalling that
in the proof of Lemma \ref{le:24} we showed that
\[ \langle p_j,\dot{q}_j \rangle \forces
\sigma_j\bnep^{-1} \circ \alpha_3\sigma_j(\bnep)  =
\check{g_2} \check{g_1},  \]
we deduce that
\[ \langle p'',\dot{q}'' \rangle \forces \check{g_1}\check{g_2}
= \check{g_2}\check{g_1}.  \]
But the equation $g_1g_2 = g_2g_1$ is about the ground model, and
hence it is true.  $\Box$
\vs

\begin{corollary}  \label{co:26}
If $(f,g_1)$ and $(f,g_2)$ are in $s$ 
then $g_1g_2^{-1}$ is in $\langle s^- \rangle$.
\end{corollary}
\vs

\proof.  There is some $g' \in \Aut (B)$ such that $(f^{-1},g')$ is in $s$.
Then by the claim, $(1,g_1g')$ and $(1,g_2g')$ are in $s$ and so $g_1g'$,
$g_2g'$ are in $s^-$.  Hence the element
\[ g_1g_2^{-1} = (g_1g')(g_2g')^{-1}  \]
lies in $\langle s^- \rangle$.  $\Box$
\vs

\begin{corollary}  \label{co:27}
Suppose $g_1, \ldots, g_k$ are elements of $\Aut (B)$ such that 
$(\nu(g_i),g_i)$ is in $s$ for each $i$, and let each of $\varepsilon_1,
\ldots,\varepsilon_k$ be either $1$ or $-1$.  If
\[\nu(g_1)^{\varepsilon_1}\ldots \nu(g_k)^{\varepsilon_k} = 1\]
then
\[g_1^{\varepsilon_1}\ldots g_k^{\varepsilon_k}\in\langle s^-\rangle.\]
\end{corollary}

\proof.  We write $f_i$ for $\nu(g_i)$.  First we show the corollary
directly in the case $k = 3$.  Taking inverses, we can assume that
$\varepsilon_2 = 1$.  When $\varepsilon_1 = \varepsilon_3 = 1$, the result
is immediate from Lemma \ref{le:24}.  We consider next the case where
$\varepsilon_1 = 1$ and $\varepsilon_3 = -1$. Here we find $g$ such that
$(f^{-1},g)$ is in $s$.  Then both of
\[  g_1^1 g_2^1 g^1,  g_3^1g^1 \]
are in $s^-$ by Lemma \ref{le:24}, and so 
\[  g_1^1g_2^1g_3^{-1} = (g_1^1g_2^1g^1)(g_3^1g_1)^{-1} \]
is in $\langle s^- \rangle$.  By symmetry this also covers the
case where $\varepsilon_1 = -1$ and $\varepsilon_3 = 1$.  Finally 
when $\varepsilon_1 = \varepsilon_3 = -1$, we repeat the same moves,
noting that
\[  g_1^{-1}g_2^1 g^1  \]
is in $\langle s^- \rangle$ by the previous case.

This case also covers the case $k=2$ by adding at the end a factor $g_3^1$
where $(1,g^3)$ is in $s$.  The case $k=1$ is trivial.

We prove the remaining cases by induction on $k$, assuming $k > 3$.  Choose
$g$ so that $(g,f_{k-1}^{\varepsilon_{k-1}}f_k^{\varepsilon_k})$ is in $s$.
Then by induction hypothesis both the elements
\[  g_1^{\varepsilon_1} \ldots g_{k-2}^{\varepsilon_{k-2}} g^1   \]
and
\[  g^{-1}g_{k-1}^{\varepsilon^{k-1}}g_k^{\varepsilon_k} \]
lie in $\langle s^- \rangle$.  Hence so does their product, 
completing the proof.  $\Box$  
\vs

Consider the subgroup $\langle s^- \rangle$ of $\Aut (B)$.  By Lemma
\ref{le:25}, $\langle s^- \rangle$ is commutative.  By Lemma \ref{le:23} and
Corollary \ref{co:27} it follows that $s$ would be a weak splitting of
$\nu$, with $\langle s^- \rangle$ as $G_0$, if for each $f$ in $\Aut (A)$
there was a unique $g$ with $\langle f, g \rangle$ in $s$.  But we can make
this true by cutting down $s$.  So $\nu$ has a weak splitting, and this
concludes the proof of Theorem \ref{th:4}.

\section{Answers to questions}   \label{se:7}
The results above answer most of the problems stated in 
\cite{hosh:1}.  In that paper we showed:
\begin{quote}
\textbf{Theorem 3 of \cite{hosh:1}\ \ } If $\cons$ is a small natural
construction in a model of ZFC, then $\cons$ is uniformisable with
parameters.
\end{quote}
We asked (Problem A) whether it is possible to remove the restriction that
$\cons$ is small.  The answer is No:

\begin{theorem}  \label{th:28}
There is a transitive model of ZFC in which some $\emptyset$-represent\-able
construction is natural but not uniformisable (even with parameters).
\end{theorem}

$\proof$.  Let $N$ be the model of Theorem~\ref{th:4}.  Let $\cons$ be some
construction $\emptyset$-representable in $N$ which is not weakly natural
(such as Example 2 in section \ref{se:2}).  Then by Theorem \ref{th:4},
$\cons$ is not uniformisable.  The rigidifying construction $\cons^r$ of
section \ref{se:2} is $\emptyset$-representable, natural and not
uniformisable.  $\Box$
\vs

Problem B asked whether in Theorem 3 of \cite{hosh:1} the formula defining
$\cons$ can be chosen so that it has only the same parameters as the
formulas chosen to represent $\cons$.  The answer is No:

\begin{theorem}  \label{th:29}
There is a transitive model $N$ of ZFC with the following property:
\begin{quote}
For every set $Y$ there are a set $X$ and a unitype rigid-based (hence small
natural) $X$-representable construction that is not $X \cup
Y$-uniformisable.
\end{quote}
\end{theorem}

\proof.  Take $N$ to be the model given by Theorem~\ref{th:4}.  Let $Y$ be
any set in $N$.  If $N$ and $Y$ are not as stated above, then for every set
$X$ and every unitype rigid-based $X$-representable construction in $N$, $X$
is $X \cup Y$-uniformisable.  So the hypothesis of Theorem \ref{th:3} holds,
and by that theorem there is in $N$ a small $\{\bar{c}\}$-uniformisable
construction that is not weakly natural.  But this contradicts the choice of
$N$.  $\Box$
\vs

Problem C asked whether there are transitive models of ZFC in which every
uniformisable construction is natural.  Theorem~\ref{th:4} is the best
answer we have for this; the problem remains open.
\vs

In \cite{ho:1} one of us asked whether there can be models of ZFC in which
the algebraic closure construction on fields is not uniformisable.

\begin{theorem}  \label{th:30}
There are transitive models of ZFC in which: 
\begin{enumerate}
\item[(a)] no formula (with or without parameters) defines for
each field a specific algebraic closure for that field, and
\item[(b)] no formula (with or without parameters) defines for each
abelian group a specific injective hull of that group.
\end{enumerate}
\end{theorem}

\proof.  Let the model $N$ be as in Theorem \ref{th:4}.  In $N$ the
constructions of Examples 3 and 4 in section \ref{se:2} are not
uniformisable, since they are not weakly natural.  So these two examples
prove (a) and (b) respectively.  $\Box$
\vs

One result in \cite{ho:1} was that there is no primitive recursive set
function which takes each field to an algebraic closure of that field.  This
is an absolute result which applies to every transitive model of ZFC, and so
it is not strictly comparable with the consistency results proved above.  In
this context we note that Garvin Melles showed \cite{me:1} that there is no
``recursive set-function'' (he gives his own definition for this notion)
which finds a representative for each isomorphism type of countable
torsion-free abelian group.

\vs

\noindent Wilfrid Hodges\\
School of Mathematical Sciences\\
Queen Mary, University of London\\
Mile End Road, London E1 4NS, England\\
\texttt{w.hodges@qmw.ac.uk}
\vs

\noindent Saharon Shelah\\
Institute of Mathematics, Hebrew University\\
Jerusalem, Israel\\
\texttt{shelah@math.huji.ac.il}

\end{document}